\documentclass[12pt]{article}

\usepackage{amsmath, epsfig, cite}
\usepackage{amssymb}
\usepackage{amsfonts}
\usepackage{latexsym}
\usepackage{float}

\newtheorem{thm}{Theorem}[section]

\newtheorem{conj}[thm]{Conjecture}
\newtheorem{lem}[thm]{Lemma}

\numberwithin{equation}{section}

\newcommand{\qed}{{\hfill$\square$}\medskip}

\begin{document}

\begin{center}
{\Large\bf Supercongruences for the $(p-1)$th Ap\'ery number
}
\end{center}

\vskip 2mm \centerline{Ji-Cai Liu$^1$ and Chen Wang$^2$}
\begin{center}
{\footnotesize $^1$Department of Mathematics, Wenzhou University, Wenzhou 325035, PR China\\
{\tt jcliu2016@gmail.com  } \\[10pt]

$^2$Department of Mathematics, Nanjing University, Nanjing 210093, PR China\\
{\tt chenwjsnu@163.com  } }
\end{center}


\vskip 0.7cm \noindent{\bf Abstract.}
In this paper, we prove two conjectural supercongruences on the $(p-1)$th Ap\'ery number, which were recently proposed by Z.-H. Sun.

\vskip 3mm \noindent {\it Keywords}:
Ap\'ery numbers; supercongruences; Bernoulli numbers
\vskip 2mm
\noindent{\it MR Subject Classifications}: Primary 11A07; Secondary 11B68, 05A19
\section{Introduction}
In 1979, Ap\'ery  \cite{apery-asterisque-1979} introduced the following numbers
\begin{align*}
A_n=\sum_{k=0}^n{n\choose k}^2{n+k\choose k}^2\quad \text{and}\quad
A_n^{'}=\sum_{k=0}^n{n\choose k}^2{n+k\choose k},
\end{align*}
in his ingenious proof of the irrationality of $\zeta(2)$ and $\zeta(3)$.
These numbers are now known as Ap\'ery numbers. Since the appearance
of these numbers, some interesting arithmetic properties have been gradually discovered.
For instance, Beukers \cite{beukers-jnt-1985} showed that for primes $p\ge 5$ and $m,r\in\mathbb{N}$,
\begin{align*}
&A_{mp^r-1}\equiv A_{mp^{r-1}-1}\pmod{p^{3r}},\\[5pt]
&A_{mp^r-1}^{'}\equiv A_{mp^{r-1}-1}^{'}\pmod{p^{3r}}.
\end{align*}
In 2012, Z.-W. Sun \cite{sunzw-jnt-2012} proved that for any prime $p\ge 5$,
\begin{align*}
\sum_{k=0}^{p-1}(2k+1)A_k\equiv p+\frac{7}{6}p^4B_{p-3}\pmod{p^5}.
\end{align*}
Here the $n$th Bernoulli number $B_n$ is defined as
\begin{align*}
\frac{x}{e^x-1}=\sum_{n=0}^{\infty}B_n\frac{x^n}{n!}.
\end{align*}

In the past two decades, some interesting congruence properties for Ap\'ery numbers and similar numbers have been widely studied (see, for example, \cite{beukers-jnt-1985,ccs-injt-2010,delaygue-cm-2018,gessel-jnt-1982,gz-jnt-2012,pan-jnt-2014,sunzh-a-2018,
sunzw-b-2013,sunzw-jnt-2012,sunzw-scm-2014}).

Our interest concerns the following two conjectural supercongruences on the $(p-1)$th Ap\'ery number, which were recently proposed by Z.-H. Sun \cite[Conjecture 2.1 \& 2.2]{sunzh-a-2018}.
\begin{conj} (Z.-H. Sun)
For any prime $p\ge 5$, we have
\begin{align}
A_{p-1}&\equiv 1+\frac{2}{3}p^3B_{p-3}\pmod{p^4},\label{a1}\\[5pt]
A^{'}_{p-1}&\equiv 1+\frac{5}{3}p^3B_{p-3}\pmod{p^4}.\label{a2}
\end{align}
\end{conj}

The aim of this paper is to prove \eqref{a1} and \eqref{a2} by establishing their generalizations.
\begin{thm}\label{t1}
Let $p\ge 7$ be a prime. Then
\begin{align}
A_{p-1}\equiv 1+p^3\left(\frac{4}{3}B_{p-3}-\frac{1}{2}B_{2p-4}\right)
+\frac{1}{9}p^4B_{p-3}\pmod{p^5}.\label{a3}
\end{align}
\end{thm}

\begin{thm}\label{t2}
Let $p\ge 7$ be a prime. Then
\begin{align}
A^{'}_{p-1}\equiv 1+p^3\left(\frac{10}{3}B_{p-3}-\frac{5}{4}B_{2p-4}\right)
+\frac{5}{18}p^4B_{p-3}\pmod{p^5}.\label{a5}
\end{align}
\end{thm}

Letting $k=1$ and $b=p-3$ in the following Kummer's congruence (see \cite[page 193]{sunzh-dam-2000}):
\begin{align*}
\frac{B_{k(p-1)+b}}{k(p-1)+b}\equiv \frac{B_b}{b}\pmod{p},
\end{align*}
we arrive at
\begin{align}
B_{2p-4}\equiv \frac{4}{3}B_{p-3}\pmod{p}.\label{a4}
\end{align}
Substituting \eqref{a4} into \eqref{a3} and \eqref{a5}, we get \eqref{a1} and \eqref{a2} for primes $p\ge 7$.
It is routine to check that \eqref{a1} and \eqref{a2} also hold for $p=5$.

The rest of this paper is organized as follows. In the next section, we prove Theorem \ref{t1}. We show Theorem \ref{t2} in the last section.

\section{Proof of Theorem \ref{t1}}
Since
\begin{align}
{p-1+k\choose k}=\frac{p}{p+k}{p+k\choose k},\label{b0}
\end{align}
we have
\begin{align}
A_{p-1}=\sum_{k=0}^{p-1}\frac{p^2}{(p+k)^2}{p-1\choose k}^2{p+k\choose k}^2.\label{b1}
\end{align}
Note that
\begin{align}
{p-1\choose k}{p+k\choose k}&=\frac{(p^2-1^2)(p^2-2^2)\cdots(p^2-k^2)}{k!^2}\notag\\
&\equiv (-1)^k\left(1-p^2H_k^{(2)}\right)\pmod{p^4},\label{b2}
\end{align}
where $H_n^{(r)}$ denotes the $n$th generalized harmonic number of order $r$:
\begin{align*}
H_n^{(r)}=\sum_{k=1}^n\frac{1}{k^r},
\end{align*}
with the convention that $H_n=H_n^{(1)}$.
It follows from \eqref{b1} and \eqref{b2} that
\begin{align}
A_{p-1}&=1+\sum_{k=1}^{p-1}\frac{p^2}{(p+k)^2}{p-1\choose k}^2{p+k\choose k}^2\notag\\
&\equiv 1+p^2\sum_{k=1}^{p-1}\frac{1-2p^2H_k^{(2)}}{(p+k)^2}\pmod{p^6}.\label{b3}
\end{align}

Furthermore, we have
\begin{align}
\frac{1}{(p+k)^2}\equiv \frac{1}{k^2}-\frac{2p}{k^3}+\frac{3p^2}{k^4}\pmod{p^3}.\label{b4}
\end{align}
Substituting \eqref{b4} into \eqref{b3} gives
\begin{align}
A_{p-1}\equiv 1+p^2H_{p-1}^{(2)}-2p^3H_{p-1}^{(3)}+3p^4H_{p-1}^{(4)}-2p^4\sum_{k=1}^{p-1}\frac{H_k^{(2)}}{k^2}\pmod{p^5}.
\label{b5}
\end{align}

For $1\le k\le p-1$, we have
\begin{align*}
H_{k}^{(2)}+H_{p-k}^{(2)}&\equiv
H_{k}^{(2)}+\sum_{i=1}^{p-k}\frac{1}{(p-i)^2}\pmod{p}\notag\\
&=H_{p-1}^{(2)}+\frac{1}{k^2}\notag\\
&\equiv \frac{1}{k^2}\pmod{p},
\end{align*}
and so
\begin{align*}
\frac{H_{k}^{(2)}}{k^2}+\frac{H_{p-k}^{(2)}}{(p-k)^2}\equiv \frac{1}{k^4}\pmod{p}.
\end{align*}
It follows that
\begin{align}
\sum_{k=1}^{p-1}\frac{H_k^{(2)}}{k^2}
&=\sum_{k=1}^{\frac{p-1}{2}}
\left(\frac{H_{k}^{(2)}}{k^2}+\frac{H_{p-k}^{(2)}}{(p-k)^2}\right)\notag\\
&\equiv \sum_{k=1}^{\frac{p-1}{2}}\frac{1}{k^4}\pmod{p}.\label{new-1}
\end{align}
By \cite[page 353]{lehmer-am-1938}, we have
\begin{align}
\sum_{k=1}^{\frac{p-1}{2}}\frac{1}{k^4}
&\equiv \frac{1}{2}\left(\sum_{k=1}^{\frac{p-1}{2}}\frac{1}{k^4}
+\sum_{k=1}^{\frac{p-1}{2}}\frac{1}{(p-k)^4}\right)\pmod{p}\notag\\
&=\frac{1}{2}H_{p-1}^{(4)}\notag\\
&\equiv 0\pmod{p}.\label{new-2}
\end{align}
It follows from \eqref{new-1} and \eqref{new-2} that
\begin{align}
\sum_{k=1}^{p-1}\frac{H_k^{(2)}}{k^2}\equiv 0\pmod{p}.\label{b6}
\end{align}
By \cite[page 353]{lehmer-am-1938}, we have
\begin{align}
H_{p-1}^{(3)}&\equiv 0\pmod{p^2},\label{b7}\\
H_{p-1}^{(4)}&\equiv 0\pmod{p}.\label{b8}
\end{align}
Combining \eqref{b5} and \eqref{b6}--\eqref{b8} gives
\begin{align}
A_{p-1}\equiv 1+p^2H_{p-1}^{(2)}\pmod{p^5}.\label{b9}
\end{align}

Letting $k=2$ in \cite[Theorem 5.1, (a)]{sunzh-dam-2000} and simplifying gives
\begin{align}
H_{p-1}^{(2)}\equiv \left(\frac{4}{3}B_{p-3}-\frac{1}{2}B_{2p-4}\right)p
+\left(\frac{4}{9}B_{p-3}-\frac{1}{4}B_{2p-4}\right)p^2\pmod{p^3}.\label{b10}
\end{align}
Substituting \eqref{a4} into \eqref{b10} yields
\begin{align}
H_{p-1}^{(2)}\equiv \left(\frac{4}{3}B_{p-3}-\frac{1}{2}B_{2p-4}\right)p
+\frac{1}{9}p^2B_{p-3}\pmod{p^3}.\label{b11}
\end{align}
The proof of \eqref{a3} follows from \eqref{b9} and \eqref{b11}.

\section{Proof of Theorem \ref{t2}}
In order to prove Theorem \ref{t2}, we need the following combinatorial identities.
\begin{lem}
For any non-negative integer $n$, we have
\begin{align}
&\sum_{k=1}^{2n}\frac{(-1)^k}{k}{2n\choose k}{2n+1+k\choose k}=-2H_{2n},\label{lw-1}\\[5pt]
&\sum_{k=1}^{n}\frac{(-1)^k}{k}{n\choose k}{n+1+k\choose k}H_k
=2\left(\sum_{k=1}^n\frac{(-1)^k}{k^2}+\frac{(-1)^n}{n+1}H_n\right).\label{wz-1}
\end{align}
\end{lem}
{\noindent\it Proof.}
By Schneider's computer algebra package {\tt Sigma} (see \cite{schneider-slc-2007}), we find that the left-hand side of \eqref{lw-1} satisfies the following recurrence:
\begin{align*}
&(-1 - 2 n)\text{SUM}[n]+2 (3 + 2 n)\text{SUM}[1+n]+(-5 - 2 n)\text{SUM}[2+n]\\
&=\frac{17 + 24 n + 8 n^2}{(1 + n) (2 + n) (3 + 2 n)}.
\end{align*}
It is easy to verify that the right-hand side of \eqref{lw-1} also satisfies the recurrence above and both sides of \eqref{lw-1} are equal for $n=0,1$.

In a similar way, both sides of \eqref{wz-1} satisfy the same recurrence:
\begin{align*}
&(1 + n) (2 + n)^2 (7 + 2 n)\text{SUM}[n]+(2 + n) (7 + 2 n) (7 + 6 n + n^2) \text{SUM}[1+n]\\
&-(3 + n) (3 + 2 n) (2 + 4 n + n^2) \text{SUM}[2+n]-(3 + n)^2 (4 + n) (3 + 2 n)\text{SUM}3+n] =0.
\end{align*}
It is routine to check that both sides of \eqref{wz-1} are equal for $n=0,1,2$.
\qed

{\noindent \it Proof of \eqref{a5}.}
By \eqref{b0} and \eqref{b2}, we have
\begin{align}
A_{p-1}^{'}&=1+\sum_{k=1}^{p-1}\frac{p}{p+k}{p-1\choose k}^2{p+k\choose k}\notag\\
&\equiv1+ p\sum_{k=1}^{p-1}\frac{(-1)^k}{p+k}{p-1\choose k}\left(1-p^2H_k^{(2)}\right)\pmod{p^5}.\label{c3}
\end{align}

Letting $n=p-1$ and $x=p$ in the following partial fraction decomposition:
\begin{align*}
\sum_{k=0}^n\frac{(-1)^k}{x+k}{n\choose k}=\frac{n!}{(x)_{n+1}},
\end{align*}
we arrive at
\begin{align*}
\sum_{k=0}^{p-1}\frac{(-1)^k}{p+k}{p-1\choose k}=\frac{1}{p{2p-1\choose p-1}}.
\end{align*}
It follows from the above that
\begin{align}
p\sum_{k=1}^{p-1}\frac{(-1)^k}{p+k}{p-1\choose k}&=
p\sum_{k=0}^{p-1}\frac{(-1)^k}{p+k}{p-1\choose k}-1\notag\\
&= \frac{1}{{2p-1\choose p-1}}-1.\label{c4}
\end{align}
We need the following McIntosh's congruence (see \cite[(6)]{mestrovic-a-2011}):
\begin{align}
{2p-1\choose p-1}\equiv 1-p^2H_{p-1}^{(2)}\pmod{p^5}.\label{c5}
\end{align}
Substituting \eqref{c5} into \eqref{c4} and using the fact that $H_{p-1}^{(2)}\equiv 0\pmod{p}$, we arrive at
\begin{align}
p\sum_{k=1}^{p-1}\frac{(-1)^k}{p+k}{p-1\choose k}\equiv p^2H_{p-1}^{(2)}\pmod{p^5}.
\label{c6}
\end{align}

On the other hand, using ${p-1\choose k}\equiv (-1)^k(1-pH_k)\pmod{p^2}$ we have
\begin{align}
&p^3\sum_{k=1}^{p-1}\frac{(-1)^k}{p+k}{p-1\choose k}H_k^{(2)}\notag\\
&\equiv p^3\sum_{k=1}^{p-1}\frac{H_k^{(2)}}{p+k}-p^4\sum_{k=1}^{p-1}\frac{H_kH_k^{(2)}}{p+k}\pmod{p^5}.\label{new-4}
\end{align}
Next, we shall prove that
\begin{align}
\sum_{k=1}^{p-1}\frac{H_kH_k^{(2)}}{p+k}\equiv \sum_{k=1}^{p-1}\frac{H_kH_k^{(2)}}{k}\equiv 0\pmod{p}.\label{new-5}
\end{align}
Letting $n=p-1$ in \eqref{wz-1} and noting \eqref{b2}, we find that
\begin{align*}
\sum_{k=1}^{p-1}\frac{\left(1-p^2H_k^{(2)}\right)H_k}{k}
\equiv 2\left(\sum_{k=1}^{p-1}\frac{(-1)^k}{k^2}+\frac{1}{p}H_{p-1}\right)\pmod{p^4}.
\end{align*}
It follows that
\begin{align*}
p^2\sum_{k=1}^{p-1}\frac{H_kH_k^{(2)}}{k}
\equiv \sum_{k=1}^{p-1}\frac{H_k}{k}-2\left(\sum_{k=1}^{p-1}\frac{(-1)^k}{k^2}+\frac{1}{p}H_{p-1}\right)\pmod{p^4}.
\end{align*}
Since
\begin{align*}
&\sum_{k=1}^{p-1}\frac{H_k}{k}=\frac{1}{2}\left(H_{p-1}^2+H_{p-1}^{(2)}\right),\\[5pt]
&\sum_{k=1}^{p-1}\frac{(-1)^k}{k^2}=\frac{1}{2}H_{(p-1)/2}^{(2)}-H_{p-1}^{(2)},\\[5pt]
&H_{p-1}\equiv 0\pmod{p^2},
\end{align*}
we have
\begin{align*}
p^2\sum_{k=1}^{p-1}\frac{H_kH_k^{(2)}}{k}
\equiv -\frac{2}{p}H_{p-1}+\frac{5}{2}H_{p-1}^{(2)}-H_{(p-1)/2}^{(2)}\pmod{p^4}.
\end{align*}
From \cite[(6) \& (7)]{mestrovic-a-2011}, we deduce that
\begin{align*}
\frac{1}{p}H_{p-1}\equiv -\frac{1}{2}H_{p-1}^{(2)}\pmod{p^3}.
\end{align*}
Thus,
\begin{align}
p^2\sum_{k=1}^{p-1}\frac{H_kH_k^{(2)}}{k}
\equiv \frac{7}{2}H_{p-1}^{(2)}-H_{(p-1)/2}^{(2)}\pmod{p^3}.\label{wz-2}
\end{align}
Letting $k=2$ in \cite[Theorem 5.2, (a)]{sunzh-dam-2000} and simplifying gives
\begin{align}
H_{(p-1)/2}^{(2)}\equiv \left(\frac{14}{3}B_{p-3}-\frac{7}{4}B_{2p-4}\right)p+\left(\frac{14}{9}B_{p-3}
-\frac{7}{8}B_{2p-4}\right)p^2\pmod{p^3}.\label{wz-3}
\end{align}
Substituting \eqref{b10} and \eqref{wz-3} into \eqref{wz-2} gives
\begin{align*}
p^2\sum_{k=1}^{p-1}\frac{H_kH_k^{(2)}}{k}
\equiv0\pmod{p^3},
\end{align*}
which implies \eqref{new-5}.

Since
\begin{align*}
\frac{1}{p+k}\equiv \frac{1}{k}-\frac{p}{k^2}\pmod{p^2},
\end{align*}
by \eqref{b6} we have
\begin{align}
\sum_{k=1}^{p-1}\frac{H_k^{(2)}}{p+k}
&\equiv \sum_{k=1}^{p-1}\frac{H_k^{(2)}}{k}-p\sum_{k=1}^{p-1}\frac{H_k^{(2)}}{k^2}\notag\\
&\equiv  \sum_{k=1}^{p-1}\frac{H_k^{(2)}}{k}\pmod{p^2}.\label{new-6}
\end{align}
Combining \eqref{new-4}, \eqref{new-5} and \eqref{new-6} gives
\begin{align}
p^3\sum_{k=1}^{p-1}\frac{(-1)^k}{p+k}{p-1\choose k}H_k^{(2)}
\equiv p^3\sum_{k=1}^{p-1}\frac{H_k^{(2)}}{k}\pmod{p^5}.\label{new-7}
\end{align}
Letting $n=\frac{p-1}{2}$ in \eqref{lw-1}, we get
\begin{align*}
\sum_{k=1}^{p-1}\frac{(-1)^k}{k}{p-1\choose k}{p+k\choose k}=-2H_{p-1}.
\end{align*}
By \eqref{b2}, we have
\begin{align*}
\sum_{k=1}^{p-1}\frac{1-p^2H_k^{(2)}}{k}\equiv -2H_{p-1}\pmod{p^4},
\end{align*}
and so
\begin{align*}
p^2\sum_{k=1}^{p-1}\frac{H_k^{(2)}}{k}\equiv 3H_{p-1}\pmod{p^4},
\end{align*}
which implies that
\begin{align}
\sum_{k=1}^{p-1}\frac{H_k^{(2)}}{k}\equiv \frac{3}{p^2}H_{p-1}\pmod{p^2}.\label{lw-2}
\end{align}
Substituting \eqref{lw-2} into \eqref{new-7} gives
\begin{align}
p^3\sum_{k=1}^{p-1}\frac{(-1)^k}{p+k}{p-1\choose k}H_k^{(2)}\equiv 3pH_{p-1}\pmod{p^5}.\label{new-9}
\end{align}
From \cite[(6) \& (7)]{mestrovic-a-2011}, we deduce that
\begin{align}
pH_{p-1}\equiv -\frac{p^2}{2}H_{p-1}^{(2)}\pmod{p^5}.\label{lw-3}
\end{align}

Finally combining \eqref{c3}, \eqref{c6}, \eqref{new-9} and \eqref{lw-3} gives
\begin{align}
A_{p-1}^{'}\equiv 1+\frac{5}{2}p^2H_{p-1}^{(2)}\pmod{p^5}.\label{new-10}
\end{align}
Then the proof of \eqref{a5} follows from \eqref{b11} and \eqref{new-10}.
\qed

\vskip 5mm \noindent{\bf Acknowledgments.}
The authors are grateful to Professor Zhi-Hong Sun for bringing the reference \cite{sunzh-dam-2000} to their attention and helpful conversations.

\end{document}